\documentclass[ECP]{ejpecp} 

\SHORTTITLE{Sharp Decoupling Inequalities for Second Moments and Variances of Sums}

\TITLE{Sharp Decoupling Inequalities for the Variances and Second Moments of Sums of Dependent Random Variables} 

\AUTHORS{%
  Victor H. de~la~Pe{\~n}a\footnote{Columbia University, United States of America.
    \EMAIL{vhd1@columbia.edu}}
  \and 
  Heyuan Yao\footnotemark[1]\footnote{Current affiliation: Northwestern University, United States of America. \BEMAIL{heyunayao2029@u.northwestern.edu}}
  \and 
  Demissie Alemayehu\footnotemark[1]\footnote{Pfizer and Columbia University, United States of America. \EMAIL{da15@columbia.edu}}
  \and Victor Keegan de~la~Pe{\~n}a\footnotemark[1]\footnote{Infremacy LLC, United States of America. \EMAIL{vkd2107@columbia.edu}}
    }


\AMSSUBJ{60E15; 60G42; 60G50; 62L15} 
\AMSSUBJSECONDARY{15A63} 

\SUBMITTED{July 22, 2025} 
\ACCEPTED{June 24, 2026} 

\VOLUME{0}
\YEAR{2026}
\PAPERNUM{0}
\DOI{10.1214/YY-TN}

\ABSTRACT{
We consider the problem of constructing sharp upper and lower bounds for the second moment of sums of arbitrary, square-integrable, dependent random variables, $\sumn d_i$, using recent results in the theory of decoupling inequalities. Specifically, we give a new proof for the complete decoupling inequality, which provides a lower bound for the sum of dependent square-integrable nonnegative random variables
\[
\frac{1}{2} \E \lpar \sumn  z_i \rpar^2 \leq \E \lpar \sumn  d_i \rpar^2,
\]
where $z_i \stackrel{\mathcal{L}}{=} d_i$ for all $i\leq n$ and $z_i$'s are mutually independent. We further develop sharp tangent decoupling inequalities for the variance and second moment of arbitrary square-integrable sequences $\{d_i\}_{1\leq i \leq n}$, without requiring nonnegativity:
\[\var \lpar \sumn d_i\rpar \leq 2 \var \lpar \sumn e_i\rpar,\]
and 
\[\E \lpar \sumn d_i\rpar^2 \leq 2 \E \lpar \sumn e_i\rpar^2 - \lbr\E \lpar \sumn e_i\rpar \rbr^2,\]
where $\{e_i\}$ is the decoupled sequences of $\{d_i\}$. These results sharpen and generalize previous bounds for randomly stopped sums and martingale-type sequences. Applications include refined Chebyshev and Paley-Zygmund-type inequalities under minimal moment assumptions. The inequalities presented here offer new tools for probabilistic analysis in both classical and modern settings involving dependent structures.
}

\newcommand{\R}{\mathbb{R}}
\newcommand{\F}{\mathcal{F}}
\newcommand{\G}{\mathcal{G}}
\newcommand{\II}{\mathbb{I}}

\newcommand{\lpar}{\left(}
\newcommand{\rpar}{\right)}
\newcommand{\lbr}{\left[}
\newcommand{\rbr}{\right]}
\newcommand{\pro}{\mathbb{P}}
\newcommand{\E}{\mathbb{E}}
\newcommand{\N}{\mathbb{N}}

\newcommand{\tX}{\tilde{X}}
\newcommand{\X}{\mathcal{X}}
\newcommand{\one}{\mathds{1}}

\newcommand{\sumn}{\sum\limits_{i=1}^n}
\newcommand{\var}{\mathbb{V}ar}

\begin{document}
\allowdisplaybreaks[1]


\section{Introduction}

\subsection{Background of Decoupling Inequalities}
As an important tool in stochastic analysis, decoupling has been studied since the 1980s, which aims to compare two random processes, where one has a weaker dependence structure. Such a comparison involves decoupling inequalities that bound the quantities of the heavier dependent structure we are interested in, such as moments and tail probabilities, via the corresponding quantities of the weakly dependent process. The application of the decoupling technique includes martingale transforms (\cite{burkholder1983geometric}), multilinear forms (\cite{mcconnell1986decoupling}), U-processes (\cite{de1992decoupling}), and stochastic integrations (\cite{kallenberg1989multiple,kwapien1991semimartingale,mcconnell1989decoupling}). For a comprehensive review of decoupling tools, we refer our readers to de la Pe\~na and Gin\'e \cite{de2012decoupling} (2012).

In this paper, we consider decoupling the sums $\sumn d_i$ of dependent random variables $\{d_i\}$. This form involves many processes and statistics one may be interested in, including U-statistics and randomly stopped processes. Specifically, two techniques of decoupling are introduced. The most direct way is the complete decoupling. As its name suggests, a sequence of mutually independent variables $\{z_i\}$ is constructed, such that $z_i$ has the same distribution as $d_i$ for all $i$. We then study if there is a constant $c$, which does not depend on any information of $\{d_i\}$ and $\{z_i\}$, such that 
\[ c\,\E  \lpar  \sumn z_i\rpar^2 \leq\E  \lpar  \sumn d_i\rpar^2.\]
Note that the lower bound does not exist in the general case. When $(d_1, d_2) $ takes values $(1,-1)$ and $(-1,1)$ with probability $1/2$, we have $\E(d_1+d_2)^2 = 0$ while $\E(e_1+e_2)^2 = 2$.  However, when we further assume that $d_i$'s are non-negative, we are able to find a lower bound given by the moment of $\sumn z_i$. The general result for any $p$-th moment was first proposed in \cite{de1990bounds} and improved and applied in \cite{makarychev2018solving}. To complete this article, we provide a new proof for the second moment case where the best constant $c$ is $\frac{1}{2}$. 
 
However, it is not possible to find a satisfactory upper bound $ C \,\E  \lpar  \sumn z_i\rpar$ for $ \E  \lpar  \sumn d_i\rpar$. Consider $d_1 = ...= d_n \stackrel{\mathcal{L}}{=} \text{Bernoulli}(p)$, for some $p\in (0,1)$, then $\E(\sumn d_i)^2 = n^2 p$ while $\E(\sumn z_i)^2 = np (np+1-p)$. By letting $n \rightarrow \infty$ and $p = o(1/n)$, $\E(\sumn d_i)^2/ \E(\sumn z_i)^2 \rightarrow \infty$. Therefore, the dependence structure needs to be modified. With the development of martingale theory, tangent decoupling becomes useful. By obtaining a decoupled version $\{e_i\}$ of $\{d_i\}$ (see Definition \ref{def4 decoupled}), we can bound the moment of $\sumn d_i$ using the moment of $\sumn e_i$. One of the earliest works was given by de la Pe\~na and Govindarajulu (1992 \cite{de1992note}), who provided sharp inequalities for the second moment of randomly stopped sums of independent variables. This result was later enhanced by Hitczenko (1994 \cite{hitczenko1994sharp}), who proved a sharp bound for the $\mathbb{L}^p$-norm of $\sumn d_i$, such that 
\begin{equation}
   \label{equ: hitczenko1994sharp}
   \lVert \sumn d_i \rVert_p^p \leq2^{p-1} \lVert \sumn e_i \rVert_p^p,
\end{equation}
where $1\leq p < \infty$ and $d_i$'s are assumed nonnegative, one application of which can be found in Feng and Yang (2019 \cite{feng2019adaptive}). In 1994, de la Pe\~na \cite{de1994bound} provided a bound for the moment generating function (M.G.F.) of $S:= \sumn d_i$ using $S':= \sumn e_i$. Given both sums being light-tailed, the following bound for all finite $\lambda$ holds:
\begin{equation}
   \label{equ: de1994bound}
   M_{S}(\lambda)\leq  \lpar  M_{S'}(2\lambda) \rpar^{\frac{1}{2}}, 
\end{equation}
where $M_X(\lambda)$ denotes the moment generating function of $X$.
This inequality led to a general class of exponential inequalities for martingales and ratios (de la Pe\~na \cite{de1999general}), which further found various applications for constructing confidence regions (e.g., \cite{rakhlin2011making,nardi2011autoregressive,wang2018stochastic,ramdas2020admissible}). The above three results provided a conjecture, that it is possible to combine \cite{de1992decoupling} and \cite{hitczenko1994sharp}, and construct a tangent decoupling inequality having a similar form to \eqref{equ: hitczenko1994sharp} and working not only for nonnegative $\{d_i\}$. Moreover, Hitczenko \cite{MR1258886} extended his result \cite{hitczenko1994sharp} to the case of real-valued variables, though the bound is not sharp. In this work, we will show that when $p=2$, a sharp upper bound for $\E (\sumn d_i)^2$ can be built, such that $\E (\sumn d_i)^2 \leq C \E (\sumn e_i)^2$, where $C=2$ is the best universal constant.

In section \ref{sec: complete decoupling}, we provide a new proof for the sharp lower bound of $\E (\sumn d_i)^2$ via complete decoupling. In Section \ref{sec: Main results}, we will demonstrate the main result, Theorem \ref{thm: decoupling inequality for 2nd moment}, which extends Hitczenko's work \cite{hitczenko1994sharp} for $\mathbb{L}^2$-norm such that the nonnegative summand assumption can be removed. Theorems \ref{thm: decoupling inequality for 2nd moment}, \ref{cor: Decoupling Inequality for the Variance} and \ref{cor: Refined decoupling inequality for 2nd moment}, provide sharp upper bounds for the variances and $\mathbb{L}^2$-norms of $\sumn d_i$, respectively. These inequalities also generalize the upper bound in \cite{de1992note} beyond randomly stopped sums of independent variables. In Section \ref{sec: Applications}, we will use some examples to show the applications of the main results. In the rest of this section, we introduce the framework of tangent decoupling and some notation.

\subsection{Definitions and constructions of decoupled sequences}
In this subsection, we briefly provide the framework of tangent decoupling. From now on, we use the notations $d_i$'s and $e_i$'s to represent the possibly dependent random variables and their corresponding decoupled random variables (see Definitions \ref{def1 tangent}, \ref{def3 conditionally independent} and \ref{def4 decoupled}). We denote $\mathcal{L}(X|\F)$ as the conditional law (distribution) of a random variable $X$ on a $\sigma$-algebra $\F$.
\begin{definition}
    \label{def1 tangent}
    Let $\{d_i\}$, $\{e_i\}$ be two sequences of random variables adapted to an increasing sequence of $\sigma$-algebras $\{\F_i\}$. Then $\{d_i\}$ is said to be \textbf{tangent} to $\{e_i\}$ with respect to $\{\F_i\}$ if for all $i$, $\mathcal{L}(d_i| \F_{i-1})$ = $\mathcal{L}(e_i|\F_{i-1})$, i.e., the conditional distributions of $d_i$ given $\F_{i-1}$ and $e_i$ given $\F_{i-1}$ are the same.
\end{definition}

\begin{definition}
    \label{def3 conditionally independent}
    A sequence $\{e_i\}$ of random variables adapted to an increasing sequence of $\sigma$-algebra $\{\F_i\}$ contained in $\F$ is said to be \textbf{conditionally independent (CI)} if there exists a $\sigma$-algebra $\mathcal{G}$ contained in $\F$ such that $\{e_i\}$ is conditionally independent given $\mathcal{G}$ and $\mathcal{L}(e_i|\F_{i-1}) =\mathcal{L}(e_i|\mathcal{G})$.
\end{definition}

\begin{definition}
    \label{def4 decoupled}
    Let $\{d_i\}$ be an arbitrary sequence of random variables, then a \textbf{ CI} sequence $\{e_i\}$ which is also tangent to $\{d_i\}$ will be called a \textbf{decoupled} version of $\{d_i\}$.
\end{definition}
The following proposition provided by Kwapien and Woiczynski (1989, \cite{MR1035249})
guarantees that we can always construct a decoupled sequence, for any possibly dependent sequence $\{d_i\}$.
\begin{proposition}
    For any sequence of random variables $\{d_i\}$, one can find a decoupled sequence $\{e_i\}$ (on a possibly enlarged probability space) which is tangent to the original sequence and in addition conditionally independent given a master $\sigma$-algebra $\mathcal{G}$. Frequently $\mathcal{G} = \sigma(\{d_i\})$.
\end{proposition}

More specifically for  any sequence $\{d_i\}$, a decoupled  tangent sequence can be obtained via the following steps:
\begin{itemize}
    \item First, we take $d_1$ and $e_1$ to be two independent copies of the same random mechanism. 
    \item With ($d_1,..,d_{i-1}$), the i-th pair of variables $d_i$ and $e_i$ comes from i.i.d. copies of the same random mechanism given $\F_{i-1} = \sigma (d_1,...,d_{i-1})$. 
\end{itemize}
And $e_i$’s are conditionally independent w.r.t. $\G = \F_n$.

An example of constructing a tangent sequence can be found in section \ref{sec: Applications}. Here we provide another example using (generalized) U-statistics. A generalized U-statistic $U_n$ of $n$ independent random variable $X_1,...,X_n \in \X$ with a class of functions $f_{i,j}: \X\times \X \rightarrow \R$, $1\leq i<j\leq n$, is defined as   
\begin{equation}
    \label{equ: generalized U stats}
    U_n= \sum\limits_{1\leq i<j\leq n}  f_{i,j}(X_i,X_j).
\end{equation}
When $f_{i,j}(\cdot,\cdot) \equiv f(\cdot,\cdot) := \frac{2}{n(n-1)} h(\cdot,\cdot)$ and $X_i$'s are i.i.d., $U_n$ becomes the U-statistic with kernel $f$. A classical way to handle the variance is Hoeffding decomposition (see e.g., Chapters 11 \& 12 \cite{MR1652247}), while decoupling considers the natural filtration $\mathbb F = \{\F_{i}\}_{1\leq i\leq n}$, where $\F_i = \sigma(\{X_1,...,X_i\})$. With $d_j : = \sum\limits_{i=1}^{j-1} f_{i,j}(X_i, X_j)$, we rewrite $U_n : = \sum\limits_{j=2}^n d_j$. The decoupled sequence $e_j: = \sum\limits_{i=1}^{j-1} f_{i,j}(X_i,\tX_j)  $ is then constructed, with $\tX_i$ the independent copy of $X_i$ for each $i$. As an instance of generalized U-statistics, the quadratic form $Q_n = \sum\limits_{1\leq i<j\leq n}  a_{ij} X_i X_j$ has the corresponding decoupled sum $\sum\limits_{j=2}^n \lpar \sum\limits_{i=1}^{j-1} a_{ij} X_i \rpar\tilde{X}_j$. With the decoupling inequalities in Section \ref{sec: Main results}, we can use the variance (or the second moment) of the decoupled sum to bound the variance (or the second moment) of $U_n$. Similar to the technique in section \ref{sec: Applications}, we can also derive an upper bound for the variance (or the second moment) of $U_\tau$, for some $\mathbb F$-stopping time $\tau$.

\section{Lower bound via complete decoupling}
\label{sec: complete decoupling}

We now provide a new proof for the following theorem, which coincides with the result in \cite{de1990bounds} and \cite{makarychev2018solving}.
\begin{theorem}
    \label{thm: complete decoupling}
    For any nonnegative dependent square-integrable  variables $d_1,...d_n$, and the mutually independent  variables $z_1,...,z_n$ such that $z_i \stackrel{\mathcal{L}}{=} d_i$ for all $i=1,...,n$, we have the following sharp inequality
    \begin{equation}
        \label{equ: complete decoupling}
        \frac{1}{2}\E\lpar\sumn z_i\rpar^2 \leq \E \lpar \sumn d_i\rpar^2.
    \end{equation}
\end{theorem}
\begin{proof}
     We notice that $\E (\sumn z_i) ^2 =  \sumn \E z_i^2 + \sum\limits_{1\leq i<j\leq n} \E [z_i z_j] $ and $\E (\sumn d_i) ^2 =  \sumn \E d_i^2 + \sum\limits_{1\leq i<j\leq n} \E [d_i d_j]$. Therefore, since all random variables we consider are nonnegative, we have 
     \[ \E (\sumn d_i) ^2 \geq \sumn \E d_i^2 =\sumn \E z_i^2, \]
    and by Jensen's inequality, 
    \[ \E (\sumn d_i) ^2  \geq (\sumn \E d_i)^2 =(\sumn \E z_i)^2 \geq \sum\limits_{1\leq i<j\leq n} \E [z_i]\E[z_j]  =\sum\limits_{1\leq i<j\leq n}  \E [z_i z_j] . \]
    By adding the above two inequalities, we claim \eqref{equ: complete decoupling}. 
    
    The sharpness of \eqref{equ: complete decoupling} can be seen in the following example. Consider $(d_1,...,d_n) = \boldsymbol{u}_k$ with probability $1/n$, for $k=1,...,n$, where $\boldsymbol{u}_k$ denotes the $n$-dimensional unit vector with the $k$-th element being 1. Therefore, $\sumn d_i \equiv 1$, $\sumn z_i \stackrel{\mathcal{L}}{=} \text{Binomial}(n,1/n)$. Hence $\E (\sumn d_i)^2 =1$ and $\E (\sumn e_i)^2 =2-1/n$. This suggests that the best constant for this lower bound working for all $n$ is $1/2$. 
\end{proof}

An application of Theorem \ref{thm: complete decoupling} is to construct a lower bound for the variance of the Horvitz–Thompson estimator \cite{MR53460}. Given an observed nonnegative sample of $n$ variables $\{y_i:i\in \mathcal S\}$, where $\mathcal S \subset \{1,...,N\}$ denotes a subset of $N$ distinct strata, we assume the population mean of  $y_1,...,y_N$ is $\mu$, and let $I_i:= \one_{\{i \in \mathcal S\}}$ denote the random variable that represents whether $i$ is present in the sample $\mathcal S$. We let $\pi_i:= \E [I_i] = \pro (i\in \mathcal S)$ denote the first-order inclusion probability and $\pi_{ij}:= \E [I_i I_j ] = \pro (i,j\in \mathcal S)$ denote the second-order inclusion probability. We assume that $\{\pi_i\}_{i=1}^N$ are known, and define the Horvitz–Thompson estimator as \[\hat \mu_{HT}: =\frac{1}{N} \sum_{i\in \mathcal S} \frac{y_i}{\pi_i} = \frac{1}{N} \sum\limits_{i=1}^N \lpar \frac{y_i}{\pi_i} I_i \rpar,\] which is an unbiased estimator of $\mu$. 

The variance of $\hat \mu_{HT}$, however, depends on the second-order inclusion probabilities $\pi_{ij} = P(i,j \in S)$, which are often intractable under complex multi-stage designs.  However, a lower bound for $\var (\hat \mu_{HT})$ can be built via our inequality \eqref{equ: complete decoupling}. A decoupled sequence of $\{ \frac{y_i}{\pi_i} I_i\}_{i=1}^N$ is given by $\{ \frac{y_i}{\pi_i} B_i\}_{i=1}^N$, where $B_i\sim \text{Bernoulli}(\pi_i)$ are independent random variables. We then obtain the decoupled Horvitz–Thompson estimator $ \tilde \mu_{HT}: \frac{1}{N} \sum\limits_{i=1}^N \lpar \frac{y_i}{\pi_i} B_i \rpar$, with mean $\mu$ and variance $\var(\tilde \mu_{HT}) = \frac{1}{N} \sum\limits_{i=1}^N  \frac{y_i^2}{\pi_i}(1-\pi_i)$. Our second moment inequality \eqref{equ: complete decoupling} then implies the variance bound, such that 
\begin{equation}
    \label{equ: LB for HT}
    \var (\hat \mu_{HT}) \geq \frac{1}{2} \var(\tilde \mu_{HT}) - \frac{1}{2}\mu^2.
\end{equation}
The bound \eqref{equ: LB for HT} is informative when the coefficient of variation $ \sqrt{\var(\tilde \mu_{HT})}/\mu>1$. This condition holds naturally in rare disease epidemiology, where the prevalence, $\mu$, is small and the inclusion probabilities $\pi_i$ for affected individuals are tiny under complex designs so that the terms $y_i^2/\pi_i$ dominate.

The sharpness of the constant $1/2$ has concrete consequences.
Makarychev and Sviridenko \cite{makarychev2018solving} showed that the $q$-th moment generalization of \eqref{equ: complete decoupling}, with sharp constant $A_q = E[P^q]$ for $P\sim\operatorname{Poisson}(1)$ established in \cite{de1990bounds}, determines the tight integrality gap of an Linear Programming (LP) relaxation for optimization problems with diseconomies of scale. In their
framework, randomized rounding produces loads $L_j = \sum_i d_{ij} R_i$ with independent Bernoulli random variables $R_i$--a Horvitz--Thompson-type sum, while the LP encodes a dependent coupling with the same marginals. The approximation ratio is exactly the decoupling constant: For $q=2$, this is $A_2 = 2$ from Theorem \ref{thm: complete decoupling}. Sharpness means the algorithm achieves this factor and no rounding based on the same relaxation can do better, yielding optimal constant-factor approximations for energy-efficient routing, load balancing, and degree-balanced spanning trees \cite{makarychev2018solving}.

\section{Upper bound for variance and second moment via tangent decoupling}\label{sec: Main results}

In this section, we present the main results of the decoupling inequalities for variance and the second moment of the sums. We assume that all random variables we consider are square-integrable.

\begin{theorem}[Moment bounds for sums of arbitrary random variables]
    \label{thm: decoupling inequality for 2nd moment}

    Consider a sequence of possibly dependent random variables $d_1,..,d_n$, and a decoupled sequence $e_1,..,e_n$, given a $\sigma$-algebra $\G$. We have the following decoupling inequality for the second moment
        \begin{equation}
            \label{equ: Decoupling Inequality for the Second Moment}
            \E \lpar \sumn d_i\rpar^2 \leq 2 \E \lpar \sumn e_i\rpar^2.
        \end{equation}
\end{theorem}

\begin{proof}
    We first use martingale-compensator decomposition, such that $\sumn d_i$ can be written as 
    \begin{equation}
        \label{equ: martingale-compensator decomposition}
        \sumn d_i = \sumn \lpar d_i-\E[d_i|\F_{i-1}]\rpar  +\sumn \E[d_i|\F_{i-1}].
    \end{equation}
    By applying the triangular inequality $(a+b)^2 \leq 2a^2+2b^2$, we upper-bound the second moment of $\sumn d_i$, such that
    \begin{align}
        \E \lpar \sumn d_i\rpar^2 & \leq 2 \E \lbr \sumn \lpar d_i-\E[d_i|\F_{i-1}]\rpar \rbr^2  + 2 \E \lbr \sumn \E[d_i|\F_{i-1}]\rbr^2 \label{equ: after tri ineq}.
    \end{align}
    The first term in \eqref{equ: after tri ineq} is the second moment of the sum of martingale differences, which can be written such that 
    \begin{equation}
        \label{equ: martingale difference}
        \E \lbr \sumn \lpar d_i-\E[d_i|\F_{i-1}]\rpar \rbr^2  =  \sumn  \E d^2_i - \E \lpar \sumn \E^2[d_i|\F_{i-1}]\rpar.
    \end{equation}

    The upper bound is then given by 
    \begin{equation}
    \label{equ: upper bound d_i}
        2 \sumn  \E d^2_i - 2\E \lpar \sumn \E^2[d_i|\F_{i-1}]\rpar  + 2 \E \lbr \sumn \E[d_i|\F_{i-1}]\rbr^2.
    \end{equation}
    By Definition \ref{def1 tangent}, we note that $d_i \stackrel{\mathcal L}{=} e_i$ and  $\E(d_i|\F_{i-1}) \stackrel{\pro\text{-a.e.}}{=} \E(e_i|\F_{i-1})$. Therefore, the upper bound \eqref{equ: upper bound d_i} can be equivalently represented by 
    \begin{equation}
        \label{equ: upper bound e_i}
        2 \sumn  \E e^2_i - 2\E \lpar \sumn \E^2[e_i|\F_{i-1}]\rpar  + 2 \E \lbr \sumn \E[e_i|\F_{i-1}]\rbr^2 = 2 \E \lpar \sumn e_i\rpar^2.
    \end{equation}
    And we conclude Theorem \ref{thm: decoupling inequality for 2nd moment}.

\end{proof}

\begin{remark}
\label{rmk: Lp bound}
    We can also use the martingale-compensator decomposition \eqref{equ: martingale-compensator decomposition} and Minkowski's inequality to obtain the following $\mathbb L^p$-norm bound, such that for any $p>1$, the following bound holds:
    \begin{equation}
        \label{equ: Lp bound}
        \|\sumn d_i \|_p \leq \| \sumn (d_i - \E[\sum e_i|\G]) \|_p + \|\E[\sumn e_i|\G]\|_p.
    \end{equation}
    We will compare this inequality with our decoupling inequality in Section \ref{sec: Applications}.
\end{remark}

A direct consequence induced by Theorem \ref{thm: decoupling inequality for 2nd moment} is the decoupling inequality for the variance.

\begin{theorem}[Decoupling inequality for the variance]
    \label{cor: Decoupling Inequality for the Variance}
    With the same setting as Theorem \ref{thm: decoupling inequality for 2nd moment}, we have that
    \begin{equation}
        \label{equ: Decoupling Inequality for the Variance}
        \var \lpar \sumn d_i\rpar \leq 2 \var \lpar \sumn e_i\rpar.
    \end{equation}
\end{theorem}

\begin{proof}
    For all $1\leq i\leq n$, when we center the random variable by letting $d'_i = d_i-\E(d_i)$, the corresponding tangent decoupled variable becomes $e'_i = e_i-\E(e_i)$. Therefore, the inequality for $\mathbb{L}^2$-norms can be extended to the inequality for variances.
\end{proof}

With Theorem \ref{cor: Decoupling Inequality for the Variance}, we can then build a refined inequality for Theorem \ref{thm: decoupling inequality for 2nd moment}, by adding $\E[(\sumn d_i)]^2$ and $\E[(\sumn e_i)]^2$ to the two sides respectively.

\begin{theorem}[Decoupling inequality for the second moment]
    \label{cor: Refined decoupling inequality for 2nd moment}
    With the same setting as Theorem \ref{thm: decoupling inequality for 2nd moment}, we have that 
    \begin{equation}
        \label{equ: Refined decoupling inequality for 2nd moment}
        \E \lpar \sumn d_i\rpar^2 \leq 2 \E \lpar \sumn e_i\rpar^2 -\lbr  \E\lpar \sumn e_i\rpar \rbr^2.
    \end{equation}
\end{theorem}

\begin{remark}[Sharpness of the three decoupling inequalities]
    We note that inequalities \eqref{equ: Decoupling Inequality for the Second Moment}, \eqref{equ: Decoupling Inequality for the Variance}, and \eqref{equ: Refined decoupling inequality for 2nd moment} are all sharp, where we have the equality holds the following extreme case: Consider $\F_0 =\{\emptyset ,\Omega\}$, $\F_1 = \F_2 = \sigma \{d_1, e_1\}$ for some centered, nontrivial, and square-integrable $d_1$ and its i.i.d. copy $e_1$. We generate $d_2$ by letting $d_2 = d_1$ $\pro$-a.e., and the $e_2$ is then forced to be $d_1$. In this case the variance of $d_1+d_2$ is $\var (2d_1) = 4\var (d_1)$, while the variance of $e_1+e_2$ is $\var (e_1+d_1) = 2\var (d_1)$.
    The above two variances are identical to the corresponding second moments, because the sums are centered. The sharpness can also be verified using Example 1 in \cite{de1992note}.
\end{remark}

Once we have the information for the second moment or the variance of $ S_n = \sumn d_i$, we have Chebyshev's inequality to bound its deviation probability. In addition, when $S_n$ is nonnegative, the lower bound for $\pro(S_n>\theta \,\E(S_n))$ can be obtained by using Paley-Zygmund inequality. However, when the second moment or the variance of $S_n$ has a complicated form or cannot be obtained, we may use the above inequalities to obtain a weaker bound with the moment information of $S'_n:= \sumn e_i$. The following theorem is the direct consequence of Theorems \ref{cor: Decoupling Inequality for the Variance} and \ref{cor: Refined decoupling inequality for 2nd moment}.

\begin{theorem}
    \label{cor: Chebyshev inequ}
    With the same setting as Theorem \ref{thm: decoupling inequality for 2nd moment}, we have a Chebyshev-type inequality for any $t>0$,
    \begin{equation}
        \label{equ: Chebyshev ineq}
        \pro\lpar |S_n - \E S_n| >t \rpar \leq \frac{2\var(S'_n)}{t^2}.
    \end{equation}
    If in addition $S_n \geq 0$ a.s.,  we have a Paley-Zygmund-type inequality for any $0<\theta<1$,
    \begin{equation}
        \label{equ: Paley-Zygmund ineq}
        \pro(S_n > \theta\, \E S_n) \geq (1-\theta)^2 \, \frac{(\E S_n)^2}{ 2 \E {S'}_n^2 - (\E S_n)^2 } = (1-\theta)^2 \, \frac{1}{ \frac{2 \E {S'}_n^2}{(\E S_n)^2} - 1 }.
    \end{equation}
\end{theorem}

\section{Bounding the second moment of randomly stopped sums}
\label{sec: Applications}
The framework of tangent decoupling may not be intuitive to understand. Therefore, we end this article with the randomly stopped processes to show how to use tangent decoupling and the application of the inequality \eqref{equ: Refined decoupling inequality for 2nd moment}.

The moments of a randomly stopped sum with martingale property can be dated back to Wald \cite{wald1944cumulative} (1944). 
In particular, Chow et al (1965, \cite{chow1965moments}) showed that if $S_n$ is the summation of the i.i.d. centered, square-integrable random variables, $X_1,...,X_n$, and if the $\mathbb{F}$-stopping time $\tau$ satisfies $\E \tau <\infty$, then 
\[\E S^2_\tau = \sigma^2 \E \tau,\]
where $\mathbb{F} = \{\F_n\}$ is the filtration such that $\F_n := \sigma(X_1,...,X_n)$. However, to deal with the second moment of the randomly stopping sum with independent non-centered random variables $S_n = \sumn X_i$, the martingale tools may not succeed. In this case, the decoupling inequality can be of help. 

For the sake of simplicity, we assumed that $X_1, X_2,...$ are i.i.d.  $\R$-valued random variables with mean $\mu \neq 0$ and variance $\sigma^2$. Suppose the $\mathbb{F}$-stopping time $\tau$ is square-integrable. Note that one can drop the assumption of identical distribution to extend the following result. We notice that $S_\tau = \sum\limits_{i=1}^{\infty} d_i$ such that $d_i = X_i \II_{\{\tau\geq i\}}$. Since $ \II_{\{\tau\geq i\}} = 1 -  \II_{\{\tau\leq i-1\}}$ is $\F_{i-1}$-measurable and therefore independent of $X_i$, we can construct a decoupled sequence $e_i = X_i' \II_{\{\tau\geq i\}}$ of $d_i$, where $X_i'$ are the i.i.d. copy of $X_i$. And $\sum\limits_{i=1}^{\infty} e_i = \sum\limits_{i=1}^{\infty} X_i' \II_{\{\tau\geq i\}}$ is identically distributed as $ \sum\limits_{i=1}^{\infty} X_i \II_{\{\tau'\geq i\}} =: S_{\tau'}$ for $\tau'$ an i.i.d. copy of $\tau$. Based on this decoupling setting, we provide the following theorem, which improves the upper bound given in \cite{de1992note}.

\begin{theorem}
    \label{thm: application to randomly stopped process}
    Let $S_n = \sumn X_i$ be the summation of $n$ i.i.d. random variables with mean $\mu$ and variance $\sigma^2$, and assume an $\mathbb F$-stopping time $\tau$ is square-integrable. By letting $\tau'$ denote a copy of $\tau$, which is independent of $\sigma(\{X_i\}_{i=1}^n)$, we have the following inequality
    \begin{equation}
        \label{equ: L2 norm of randomly stopped sums}
        \E S^2_{\tau}\leq 2 \E S^2_{\tau'} - (\E S_{\tau'})^2,
    \end{equation}
    where the bound can be explicitly written as \eqref{equ: upper bound randomly stopped process}.
\end{theorem}

\begin{proof}
    Fatou's lemma shows that 
    \[ \E S^2_\tau = \E \lim\limits_{N\rightarrow \infty} \lpar \sum\limits_{i=1}^N X_i \one_{\{\tau \geq i\}} \rpar^2 \leq \liminf\limits_{N\rightarrow \infty}  \E \lpar \sum\limits_{i=1}^N X_i \one_{\{\tau \geq i\}} \rpar^2.\]
    Theorem \ref{cor: Refined decoupling inequality for 2nd moment} further provides the bound such that for any $N\in \N$,
    \begin{align*}
         \E &\lpar \sum\limits_{i=1}^N X_i \one_{\{\tau \geq i\}} \rpar^2   \leq 2\E \lpar \sum\limits_{i=1}^N \tX_i \one_{\{\tau \geq i\}} \rpar^2 - \lbr \E \lpar \sum\limits_{i=1}^N X_i \one_{\{\tau \geq i\}}\rpar \rbr^2  \\
         & = 4\E \lpar \sum\limits_{j=2}^N (\sum\limits_{i=1}^{j-1} X_i) X_j \II_{\{\tau'\geq j\}} \rpar + 2 \E \lpar \sum\limits_{i= 1}^N X_i^2 \II_{\{\tau'\geq i\}} \rpar - \lbr \E \lpar \sum\limits_{i=1}^N X_i \one_{\{\tau \geq i\}}\rpar \rbr^2 \\
         & = 4 \mu^2  \sum\limits_{j=2}^N (j-1) \pro(\tau \geq j) +  (2\mu^2+2\sigma^2) \sum\limits_{i=1}^N \pro(\tau\geq i)  -\mu^2\lpar  \sum\limits_{i=1}^N  \pro(\tau\geq i)\rpar^2,
    \end{align*}
    the limit of which, due to the square-integrability of $\tau$, is 
    \begin{equation}
        \label{equ: upper bound randomly stopped process}
        4 \mu^2  \sum\limits_{j=2}^\infty (j-1) \pro(\tau \geq j) +  (2\mu^2+2\sigma^2) \E \tau - \mu^2 (\E \tau)^2= 2 \mu^2 \E\tau^2 + 2\sigma^2 \E \tau - \mu^2 (\E \tau)^2.
    \end{equation}
    where the equality holds when we apply the tail sum formula.
    
    We will then show $2 \E S^2_{\tau'} - (\E S_{\tau'})^2 =$ \eqref{equ: upper bound randomly stopped process}. We note that
    \begin{align}
        \E S^2_{\tau'} & = \E \lbr \lim\limits_{N\rightarrow\infty} \lpar \sum \limits_{i=1}^N X_i \one_{\{\tau'\geq i\}}\rpar^2 \rbr \nonumber\\
        & = \E \lbr \lim\limits_{N\rightarrow\infty} \lpar \sum\limits_{i=1}^N X_i^2 \one_{\{\tau'\geq i\}}+2\sum\limits_{1\leq i<j\leq N} X_i X_j \one_{\{\tau'\geq j\}} \rpar \rbr\nonumber \\
        & =   \lim\limits_{N\rightarrow\infty}  \E\lpar \sum\limits_{i=1}^N X_i^2 \one_{\{\tau'\geq i\}}\rpar + 2\lim\limits_{N\rightarrow\infty} \E \lpar \sum\limits_{1\leq i<j\leq N} X_i X_j \one_{\{\tau'\geq j\}} \rpar \label{equ: due to DCT}\\
        & = \E X^2 \E \tau' + 2\mu^2  \sum\limits_{j =2}^\infty \sum\limits_{i=1}^{j-1} \pro(\tau' \geq j)  \label{equ: by independence of Xi Xj and tau'}\\
        & = (\mu^2+\sigma^2) \E \tau + 2\mu^2 \sum\limits_{j=2}^\infty (j-1) \pro(\tau \geq j),\nonumber 
    \end{align}
    where we obtain the first term in \eqref{equ: due to DCT} via the monotone convergence theorem, and the second term in \eqref{equ: due to DCT} via the dominated convergence theorem (since $|\sum X_i X_j \one_{\tau' \geq j}| \leq \sum |X_i| |X_j| \one_{\tau' \geq j}$ is integrable), and we use the mutual independence among $X_i, X_j$ and $\tau'$ to obtain \eqref{equ: by independence of Xi Xj and tau'}. 

    With the implicit form of $\E S^2_{\tau'}$ and that $\E S_{\tau'}= \mu \E \tau'$, we conclude the bound \eqref{equ: upper bound randomly stopped process}.
\end{proof}

\begin{remark}
    For any $p>1$, we next provide an $\mathbb L^p$-norm bound of $S_\tau$, inspired by \eqref{equ: Lp bound} in Remark \ref{rmk: Lp bound}. Recall that $d_i = X_i \one_{\{\tau \geq i\}}$, $e_i = X_i' \one_{\{\tau \geq i\}}$, we have that \[\E[e_i | \G] = \E[e_i | \F_{i-1}] = \mu \one_{\{\tau \geq i\}}.\]

    With some assumption on the integrability (e.g., $X_i,\tau \in \mathbb L^p$), \eqref{equ: Lp bound} implies that
    \begin{equation}
        \label{equ: Lp bound for randomly stopped time}
        \| S_\tau \|_p \leq \| \check S_\tau \|_p + |\mu|0 \|\tau \|_p,
    \end{equation}
    where $\check S_n = \sumn (X_i - \mu)$ denotes the mean corrected sum of $X_i$'s. Specifically, when $p=2$ and under the assumption of Theorem \ref{thm: application to randomly stopped process}, we have
    \begin{equation}
    \label{equ: Lp bound when p=2}
        \E S^2_\tau = \| S_\tau \|_2^2 \leq \mu^2 \mathbb E\tau^2 + \sigma^2 \mathbb E \tau + 2\sigma |\mu| \sqrt{\mathbb E \tau \mathbb E \tau^2}.
    \end{equation}
    To compare the bound \eqref{equ: Lp bound when p=2} with the bound \eqref{equ: upper bound randomly stopped process} in Theorem \ref{thm: application to randomly stopped process}, we notice that each bound has its own advantage, where the bound \eqref{equ: upper bound randomly stopped process} is better when $\mu \neq 0$ and 
    \begin{align*}
        & 2 \mu^2 \mathbb E\tau^2 + 2\sigma^2 \mathbb E \tau - \mu^2 (\mathbb E \tau)^2 < \mu^2 \mathbb E\tau^2 + \sigma^2 \mathbb E \tau + 2\sigma |\mu| \sqrt{\mathbb E \tau \mathbb E \tau^2}\\
        \Leftrightarrow&  \lpar  |\mu| \sqrt{\mathbb E\tau^2} - \sigma \sqrt{\mathbb E \tau} \rpar^2 < (\mu \mathbb E \tau)^2 \\
         \Leftrightarrow& \left | \sqrt{\frac{\E \tau^2}{\E^2 \tau} } -  \frac{\sigma}{|\mu|}\sqrt{\frac{1}{\E \tau} }  \right | < 1 \\
        \Leftrightarrow & \left | \sqrt{SCV_\tau +1 } -  \sqrt{\frac{SCV_X}{\E \tau}}  \right | < 1,
    \end{align*}
    where $SCV_\tau$ (resp., $SCV_X$) denotes the squared coefficient of variation $\mathbb Var(\tau)/(\mathbb E \tau)^2$ (resp., $\mathbb Var(X)/(\mathbb E X)^2$). 
\end{remark}
    

\section{Concluding Remarks}
\label{sec: concluding remarks}
This paper is concerned with the  development and application of sharp decoupling inequalities, with a focus on the second moments and variances of sums of dependent random variables.

Our findings demonstrate that precise upper bounds can be obtained under the minimal condition of square integrability. These inequalities not only deepen theoretical understanding but also yield practical tools for bounding tail probabilities via Chebyshev- and Paley-Zygmund-type inequalities in the presence of dependence.

We further highlight the applicability of these decoupling techniques to randomly stopped sums, showcasing their effectiveness in settings where traditional martingale methods may be inadequate, particularly in the analysis of non-centered or heavy-tailed variables.

By synthesizing existing literature and extending foundational inequalities, this work offers a meaningful contribution to the study of dependent random variables. Future research may build upon this foundation to explore applications in statistical inference and machine learning, where dependence among data points is often intrinsic. It may be fruitful to refer to the previous work of Hitczenko \cite{MR1258886} as a possible starting point. In addition, our mean-corrected technique might be combined with this framework to extend even-moment inequalities to bounds for general even-centered moments (e.g., kurtosis).





\providecommand{\bysame}{\leavevmode\hbox to3em{\hrulefill}\thinspace}
\providecommand{\MR}{\relax\ifhmode\unskip\space\fi MR }
\providecommand{\MRhref}[2]{%
  \href{http://www.ams.org/mathscinet-getitem?mr=#1}{#2}
}
\providecommand{\href}[2]{#2}


\begin{acks}
The first author would like to thank the support from Google DeepMind project number GT009019. In addition, the authors acknowledge the assistance of LLMs to improve grammar and readability. The authors would also like to appreciate the valuable suggestions and comments of two anonymous reviewers and the editor.
\end{acks}


\end{document}